\documentclass[11pt,a4paper]{article}
\usepackage[utf8]{inputenc}
\usepackage{geometry}
\usepackage{graphicx}
\usepackage{amsmath,amssymb}
\usepackage{verbatim} 
\usepackage{dsfont}
\usepackage{color}
\begin{document}

\begin{titlepage}
\begin{flushright}
February 21, 2021
\end{flushright}
\begin{flushright}
\end{flushright}
\vfill
\begin{center}
{\Large\bf Solution of the Basel problem\\in the framework of distribution theory}
\vfill
{\bf Andreas Aste}\\[1cm]
{\it Department of Physics, University of Basel,\\
Klingelbergstrasse 82, CH-4056 Basel, Switzerland\\
E-Mail: andreas.aste@unibas.ch}
\end{center}
\vfill
\begin{abstract}
A simple proof of Euler’s formula which states that the sum of the reciprocals of all
natural numbers squared equals $\pi^2/6$ is presented based on the distribution
theory introduced by Laurent Schwartz. Additional identities are obtained as a
byproduct of the derivation.
\end{abstract}
\vfill
{\bf{Mathematics Subject Classification MSC 2010:}} 40A25, 46F05
\vfill
\end{titlepage}

\section{Introduction}
The so-called {\emph{Basel problem}} to determine the sum
$\zeta(2) = \sum_{n=1}^\infty \frac{1}{n^2} = \pi^2/6$ was first posed in 1644 by
Pietro Mengoli, an Italian mathematician and clergyman from Bologna, and solved by the
Swiss mathematician Leonhard Euler (*1707 in Basel, $^\dagger$1783 in Saint
Petersburg) in 1735. Several ways have been found in the meantime to calculate
$\zeta(2)$ (see \cite{Overview} and references therein).
A further simple method to derive Euler's result using
the the theory of distributions and test functions which is based on elementary arguments
like translational invariance is presented in this letter.\\

Distribution theory \cite{Schwartz}, which represents a mathematical discipline in its own right,
is of fundamental significance for a rigourous treatment of quantum field theories in classical
spacetime \cite{PCT,Scharf}. It is also hoped that the stunning exercise presented in this letter serves as an
incentive for graduate students with some basic knowledge of distribution theory to study the
subject of generalized functions and their applications in theoretical physics in greater detail.

\section{Calculating $\zeta (2)$}
We consider the distribution $\Delta_0 \in \mathcal{D}' (\mathds{R})$ defined by the
formal expression
\begin{equation}
\Delta_0 (x) = \sum \limits_{n = -\infty}^\infty e^{inx} = \ldots + e^{-3ix} + e^{-2ix}
+ e^{-x} + 1 + e^{ix} + e^{2ix} + e^{3ix} +\ldots\, ,
\label{formal_def}
\end{equation}
which acts on (smooth) test functions (with compact support)
$\varphi \in \mathcal{D} (\mathds{R})$ as a linear
and, in the sense of distributions, continuous functional according to
\begin{equation}
\Delta_0 [\varphi] := \sum \limits_{n = -\infty}^\infty \int \limits_{-\infty}^\infty
e^{inx} \varphi (x) dx = \lim_{N \rightarrow \infty}
\sum \limits_{n = -N}^N \int \limits_{-\infty}^\infty e^{inx} \varphi (x) dx\, .
\label{definition_R}
\end{equation}
In fact, $\Delta_0$ is well-defined by equation (\ref{definition_R}) as a distribution
in $\mathcal{D}' (\mathds{R})$, the dual space of $\mathcal{D} ( \mathds{R})$, and
equation (\ref{definition_R}) highlights the meaning of the formal definition (\ref{formal_def})
of $\Delta_0$ as a {\emph{generalized}} function \cite{Constantinescu}. Note that a more intuitive
representation of $\Delta_0$ as an alternative infinite sum of Dirac delta distributions is motivated in the appendix.\\

By definition, $\Delta_0$ is a periodic distribution invariant under a translation
$T_{2 \pi}$, i.e. formally
\begin{equation}
(T_{2 \pi} \Delta_0 ) (x) = \Delta_0  (x+2 \pi) 
= \sum \limits_{n = -\infty}^\infty e^{in(x+2 \pi)}
= \sum \limits_{n = -\infty}^\infty e^{inx} = \Delta_0 (x) \, ,
\end{equation}
or in distributional notation
\begin{equation}
(T_{2 \pi} \Delta_0) [\varphi] = \Delta_0 [ T_{2 \pi} \varphi] \quad
\forall \varphi \in \mathcal{D}(\mathds{R}) \, , \quad \mbox{where} \, \, \,
(T_{2 \pi} \varphi) (x) = \varphi(x-2 \pi) \, ,
\end{equation}
and $\Delta_0$ is symmetric
\begin{equation}
\Delta_0 (x) = \Delta_0 (-x) \, .
\end{equation}

Now since $\Delta_0(x)$ is invariant with respect to a multiplication with $e^{ix}$, i.e.
\begin{equation}
e^{ix} \Delta_0 (x) = e^{ix} \sum \limits_{n = -\infty}^\infty e^{inx} = 
\sum \limits_{n = -\infty}^\infty e^{inx} \, ,
\end{equation}
$\Delta_0$ must vanish as a distribution on
$\mathds{R} \backslash \{2 \pi n \, | \, n \in \mathds{Z} \}$, since only for
$x=2 \pi n$ with $n \in \mathds{Z}$ one has a trivial factor $e^{ix}=1$; therefore
the distributional {\emph{support}} of $\Delta_0$ must be contained in a
corresponding discrete set
\begin{equation}
supp \, \Delta_0 \subseteq \{2 \pi n \, | \, n \in \mathds{Z} \} \, .
\end{equation}

For a moment, the following considerations are restricted to the open interval $I=(0, 2 \pi)$.
Calculating the first antisymmetric antiderivative $\Delta_1$ of $\Delta_0$ with $x \in I$
\begin{equation}
\Delta_1 (x) = \lim_{\epsilon \searrow 0} \int \limits_\epsilon^x \Delta_0 (x') dx' =
- i \sum \limits_{\substack{n = -\infty \\ n \neq 0}}^\infty \frac{e^{inx}}{n} + x
= 2 \sum \limits_{n=1}^\infty \frac{\sin (nx)}{n} + x
\label{Delta1_Fourier}
\end{equation}
with
\begin{equation}
\Delta_1 (x) = -\Delta_1 (-x) \, ,
\end{equation}
$\Delta_1$ must be {\emph{constant}} on $I$,
since its derivative $\Delta_0$ {\emph{vanishes}} there.
This also implies that the Fourier sum in equation (\ref{Delta1_Fourier})
represents a linear function on $I$. Calculating the {\emph{mean value}} $\mu_{I,1}$
of $\Delta_1$ on $I$ according to
\begin{equation}
\mu_{I,1} = \frac{1}{2 \pi} \lim_{\epsilon \searrow 0}
\int \limits_\epsilon^{2 \pi-\epsilon} \Delta_1 (x) dx \, ,
\end{equation}
the oscillatory terms $\sim e^{inx}$ in equation (\ref{Delta1_Fourier}) do not contribute
to $\mu_{I,1}$ and one is left with
\begin{equation}
\mu_{I,1} = \frac{1}{2 \pi} \int \limits_0^{2 \pi} x dx = \pi \, .
\end{equation}
Finally turning to the antiderivative of $\Delta_1$ on $I$
\begin{equation}
\Delta_2 (x) = \int \limits_0^x \Delta_1 (x') dx' =
-  \sum \limits_{n \in \mathds{Z} \backslash 0} \frac{e^{inx}}{n^2} + \frac{1}{2} x^2
= -2 \sum \limits_{n=1}^\infty \frac{\cos (nx)}{n^2} + \frac{1}{2} x^2  \, ,
\label{Delta2_Fourier}
\end{equation}
one arrives at an expression containing a series that converges absolutely to a
continuous function on $I$. However, since the distributional derivative of $\Delta_2$
is $\Delta_1$ which is constant, i.e., $\pi$ on $I$, $\Delta_2$ must be of the form
\begin{equation}
\Delta_2 (x) = \pi x + \gamma \, , \quad x \in I 
\end{equation}
with an integration constant $\gamma$. This constant can be calculated by considering
the average value of $\Delta_2$ on $I$
\begin{equation}
\mu_{I,2} = \frac{1}{2 \pi} \int \limits_0^{2 \pi} x^2 \frac{dx}{2} = \frac{2 \pi^2}{3}
=  \frac{1}{2 \pi} \int \limits_0^{2 \pi}  (\pi x +\gamma ) dx = \pi^2 + \gamma \, ,
\end{equation}
hence $\gamma = - \pi^2/3$, an finally Euler's famous result
\begin{equation}
\Delta_2 (0) = -2 \sum \limits_{n=1}^\infty \frac{1}{n^2} = \gamma = - \frac{\pi^2}{3}
\, \, \rightarrow \, \, \sum \limits_{n=1}^\infty \frac{1}{n^2} = \frac{\pi^2}{6}
\end{equation}
follows.\\

As an exercise, the reader may verify that by considering additional antiderivatives
of $\Delta_2$ like $\Delta_4$, $\Delta_6$ {\emph{et cetera}}, further values of
the Euler-Riemann zeta function like
\begin{equation}
\zeta(4) = \sum \limits_{n=1}^\infty \frac{1}{n^4} = \frac{\pi^4}{90} \, , \quad
\zeta(6) = \sum \limits_{n=1}^\infty \frac{1}{n^6} = \frac{\pi^6}{945} \, , \quad \ldots
\end{equation}
follow directly from strategy outlined above.

\appendix
\section{Explicit representation of $\Delta_0$ as an infinite sum of Dirac delta distributions}
\begin{figure}
\begin{center}
\includegraphics[width=15.0cm]{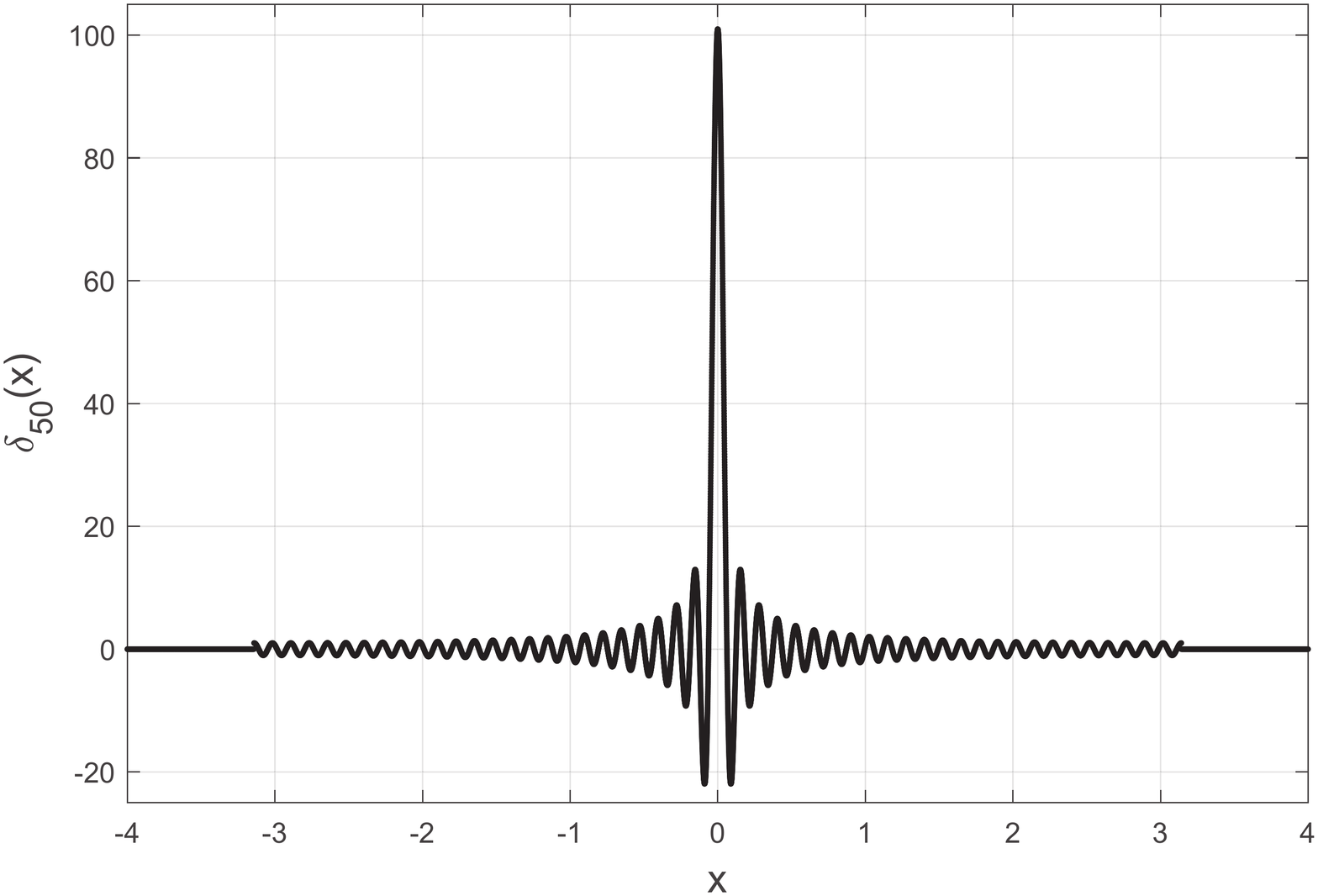}
\vskip -0.6cm
\caption{The graph of $\delta_{50}$ defined by equation (\ref{definition_delta}).}
\label{delta50}
\end{center}
\end{figure}

We consider the following sequence $\{ \delta_N \}_{N \in \mathds{N}_0} \subset
\mathcal{D}' (\mathds{R})$ 
of distributions \cite{Constantinescu} represented by the functions
\begin{equation}
\delta_N (x) = \Theta(\pi^2 - x^2) \sum \limits_{n = -N}^N e^{inx} =
    \begin{cases} 
      \sum \limits_{n = -N}^N e^{inx} & |x| < \pi \\
      0 & |x| \ge \pi 
   \end{cases}
   \label{definition_delta}
\end{equation}
with $supp \, \delta_N = [ - \pi , \pi ]$, where
\begin{equation}
\Theta (x) =
   \begin{cases} 
      1 & x >0 \\
      0 & x \le 1 
   \end{cases}
   \end{equation}
is the Heaviside function.
With $\delta_N (x) = e^{-iNx} + e^{-i (N-1)x} + \ldots + e^{-ix} + 1 + e^{ix} +
\ldots + e^{iNx}$
and
\begin{equation}
e^{ix} \delta_N (x) = \delta_N (x) + e^{i(N+1)x} - e^{-i N x}
\end{equation}
for $x \in (-\pi,\pi)$
one immediately obtains the compact representation
\begin{equation}
\delta_N (x) = \frac{e^{i (N+1)x}- e^{-iNx}}{e^{ix} -1}
= \frac{e^{i (N+1/2)x}- e^{-i(N+1/2)x}}{e^{ix/2} -e^{-ix/2}}
= \frac{\sin((N+1/2)x)}{\sin (x/2)} \, , \quad x \in (-\pi, \pi) \backslash \{ 0 \}
\label{compact_rep}
\end{equation}
and from the definition (\ref{definition_delta}) one has $\delta_N(0) = 2N + 1$ which
removes the singularity appearing at $x=0$ in the representation (\ref{compact_rep}).
Only the term $e^{i0x}=1$ in definition (\ref{definition_delta}) contributes
to the integral
\begin{equation}
\int \limits_{-\infty}^\infty \delta_N (x) dx = \int \limits_{-\pi}^\pi \delta_N (x) dx = 2 \pi \, .
\label{delta_norm}
\end{equation}
For illustrative purposes, the graph of $\delta_{50}$ is depicted in Fig. \ref{delta50}.
In fact, $\{ \delta_N \}_{N \in \mathds{N}_0} \subset
\mathcal{D}' (\mathds{R})$ is a {\emph{$\delta$-sequence}} converging to $2 \pi$ times
the Dirac delta distribution $\delta$ for $N \rightarrow \infty$.
Applying $\delta_N$ on a (smooth) test function $\varphi \in \mathcal{D} (\mathds{R})$
(with compact support) leads to
\begin{displaymath}
\delta_N [\varphi] = \int \limits_{-\infty}^\infty \delta_N (x) \varphi (x) dx
= \int \limits_{-\pi}^\pi \frac{\sin((N+1/2)x)}{\sin (x/2)} \varphi (x) dx
\end{displaymath}
\begin{equation}
= 2 \int \limits_{-\pi}^\pi \frac{\sin((N+1/2)x)}{x}
\frac{x/2}{\sin (x/2)} \varphi (x) dx
=2 \int \limits_{-\pi}^\pi \frac{\sin((N+1/2)x)}{x} \beta(x)
\frac{x/2}{\sin (x/2)} \varphi (x) dx
\label{integral_dist}
\end{equation}
where a smooth bump function $\beta \in \mathcal{D} (\mathds{R})$ with the
properties $\beta (x) = 1$ for $|x| \le \pi$ and $\beta(x) =0$ for $|x| \ge 3 \pi/2$
was introduced which does not change the integral above. Since one has
\begin{equation}
\sigma (x) =
\begin{cases} 
      \frac{1}{{\rm{si}} (x/2)} = \frac{x/2}{\sin (x/2)}  & |x| \in (0, 3 \pi/2] \\
      1 & x = 0 
   \end{cases} \, , \quad \sigma \in C^\infty ([-3 \pi / 2, 3 \pi /2]) \, ,
\end{equation}
i.e. since $\sigma$ is a smooth function on the interval $[-3 \pi / 2, 3 \pi /2]$, also 
$\tilde{\varphi}(x) = \beta(x) \sigma (x) \varphi(x)$ is smooth and has compact support:
$\tilde{\varphi} \in \mathcal{D} ( \mathds{R})$. Furthermore, $\tilde{\varphi} (0) =
\varphi(0)$ holds.\\

Now, equation (\ref{integral_dist}) becomes, with $x' = (N+1/2)x$ in the limit
$N \rightarrow \infty$ in the sense of distributions
\begin{displaymath}
\delta_N [\varphi] = \int \limits_{-\infty}^\infty \delta_N (x) \varphi (x) dx
=2 \int \limits_{-\pi}^\pi \frac{\sin((N+1/2)x)}{x} \tilde{\varphi} (x) dx
=2 \int \limits_{-(N+1/2) \pi}^{(N+1/2) \pi} \frac{\sin(x')}{x'} \tilde{\varphi} (x'/(N+1/2)) dx'
\end{displaymath}
\begin{equation}
\, \xrightarrow{N\rightarrow\infty} \,
2 \int \limits_{-\infty}^\infty \frac{\sin(x')}{x'} \tilde{\varphi} (0) dx'
= 2 \pi \tilde{\varphi}(0) = 2 \pi \varphi(0)
= 2 \pi \delta [\varphi ] \, .
\end{equation}
The normalization of the $\delta$-distribution follows from equation (\ref{delta_norm}),
i.e., as a byproduct of the derivation presented above the integral
\begin{equation}
\lim_{N \rightarrow \infty} \int \limits_{-(N+1/2)\pi}^{(N+1/2) \pi} \frac{\sin(x)}{x} dx
= \int \limits_{-\infty}^\infty \frac{\sin(x)}{x} dx = \pi
\end{equation}
is obtained.
Neglecting the cutoff in definition (\ref{definition_delta}) leads to the periodic distributional
identity
\begin{equation}
\Delta_0 (x) = \sum \limits_{n = -\infty}^\infty e^{inx}
= 2 \pi \sum \limits_{n = -\infty}^{\infty} \delta (x - 2 \pi n) \quad
\mbox{or} \quad \Delta_0 [\varphi] = 2 \pi \sum \limits_{n = -\infty}^\infty \varphi(2 \pi n) \, .
\end{equation}

One readily expresses the antisymmetric antiderivative of $\Delta_0$ by the help of the
floor function $\lfloor \cdot \rfloor$ and the ceiling function $\lceil \cdot \rceil$
\begin{equation}
\Delta_1 (x) = \pi \biggl( \biggl\lfloor \frac{x}{2 \pi} \biggr\rfloor
+ \biggl\lceil \frac{x}{2 \pi} \biggr\rceil\biggl) \, ,
\end{equation}
which simplifies to
\begin{equation}
\Delta_1 (x) = \pi \, \text{sign} (x)
\end{equation}
on the open interval $(-2 \pi , 2 \pi)$,
and the symmetric antiderivative of $\Delta_1$ is represented by the continuous function
\begin{equation}
\Delta_2 (x) = \pi x \biggl( \biggl\lfloor \frac{x}{2 \pi} \biggr\rfloor +
\biggl\lceil \frac{x}{2 \pi} \biggr\rceil \biggl) - 2 \pi^2 \biggl(
\biggl\lceil \frac{x}{2 \pi} \biggr\rceil \biggl\lfloor \frac{x}{2 \pi} \biggr\rfloor \biggl)
-\frac{\pi^2}{3} = -\sum \limits_{n \in \mathds{Z} \backslash \{ 0 \}} \frac{e^{inx}}{n^2}
+\frac{x^2}{2}\, .
\end{equation}

\end{document}